\documentclass[11pt,a4paper]{article}
\usepackage{graphicx}
\usepackage{latexsym}
\usepackage[all]{xy}
\usepackage{amsfonts}
\usepackage{amsthm}
\usepackage{amsmath}
\usepackage{amssymb}

\def\<{{\langle}}
\def\>{{\rangle}}

\def\eps{\epsilon}

\def\note#1{{}}

\def\note#1{}

\def\rend#1#2{{{\rm End}\sb{#1}(#2)}}

\def\beq{\begin{equation}}
\def\eeq{\end{equation}}

\def\ot{{\otimes}}

\newcounter{zlist}

\def\Label#1{\label{#1}\ifmmode\llap{[#1] }\else 
\marginpar{\smash{\hbox{\tiny [#1]}}}\fi}
\def\Label{\label}

\newtheorem{proposition}{Proposition}[section]

\newtheorem{theorem}[proposition]{Theorem}

\theoremstyle{definition}
\newtheorem{definition}[proposition]{Definition}

\theoremstyle{remark}
\newtheorem{remark}[proposition]{Remark}

\newcounter{c}

\newcommand{\etyk}[1]{\vspace{-7.4mm}$$\begin{equation}\Label{#1}
\addtocounter{c}{1}}
\renewcommand{\]}{\ifnum \value{c}=1 $$\else \end{equation}\fi}
\begin{document}

\title{ \bf
ON JORDAN (CO)ALGEBRAS}
\author{FLORIN F. NICHITA\\
Institute of Mathematics "Simion Stoilow" of the Romanian Academy \\ 
P.O. Box 1-764, RO-014700 Bucharest, Romania \\
E-mail: Florin.Nichita@imar.ro
}

\date{}
\maketitle
\begin{abstract}
We present new results about Jordan algebras and Jordan coalgebras, and
we discuss about their
connections with the
Yang-Baxter equations. 

\end{abstract}

{\centering\section{INTRODUCTION AND PRELIMINARIES}}

Jordan algebras emerged in the early thirties with Jordan's papers \cite{j1, j2, j3} on the algebraic formulation of quantum mechanics. Their applications
are in differential geometry, ring geometries, physics, quantum groups,
analysis, biology, etc (see \cite{RI, I}). The current paper 
 started as a poster
presented at the 11-th International Workshop on Differential Geometry and its 
Applications, in September 2013, at the Petroleum-Gas University
from Ploiesti. It contains new results on Jordan algebras and Jordan 
coalgebras, and it 
attempts to
present the general framework in which they are related to the 
quantum Yang--Baxter equation (QYBE).

Since the apparition of the QYBE
in theoretical physics (\cite{Yang}) and 
statistical mechanics (\cite{Bax1, Bax2}), many areas of mathematics,
physics and computer science have been enhanced:
knot theory, non-commutative geometry,
quantum groups,
analysis
of integrable systems, 
quantum and statistical mechanics, quantum computing,
etc (see \cite{ax}). 
Non-additive solutions of the two-parameter form of the QYBE are  
related to  the solutions of the one-parameter form of the
Yang-Baxter equation; the 
theory of integrable Hamiltonian
systems makes great use of the solutions of the one-parameter form of the
Yang-Baxter equation.


In the next section we prove a new theorem about Jordan algebras,
we explicitly define Jordan coalgebras,
and
we present a dual version of the above theorem. 
Note that Jordan coalgebras have been studied before in \cite{acm, ac},
but an explicit definition for them was never given.
Section 3 is a survey on the QYBE, and it ends
with some directions
of study related to Jordan algebras and related topics.
Our preferred bibliography on
QYBE consists of the following: \cite{LamRad:Int, Kas:Qua, HlaSno:sol, HlaKun:qua, bn, np}.

\bigskip

Throughout this paper $ k $ is a field, and all tensor products appearing are defined over $k$.
For $ V $ a $ k$-space, we denote by
$ \   \tau : V \otimes V \rightarrow V \ot V \  $ the twist map defined by $ \tau (v \ot w) = w \ot v $, and by $ I: V \rightarrow V $
the identity map of the space V.

\bigskip

{\centering\section{JORDAN ALGEBRAS AND JORDAN COALGEBRAS}}

If we consider an associative operation on a set $V$, 
then the elements of this set satisfy the Jordan identity:   

\begin{equation} \label{trei}
( x^2 y ) x = x^2 (yx) \ \ \ \ \forall x, y \in V.
\end{equation}
                                                       
If we start with an operation on a set V satisfying (\ref{trei}), 
then the elements of this set might not satisfy the associativity axiom:

\begin{equation} \label{patru}
( ab ) c = a (bc) \ \ \ \ \forall a,b,c \in V.
\end{equation}

The following question arises: 
"What conditions should we impose to an operation on a set V, 
such that if the elements of V satisfy (\ref{trei}), 
then they also satisfy (\ref{patru}) ?".

We will give an answer to this question below.

Let us recall that a Jordan algebra consists of a vector space $V$ and a 
linear map $ \theta : V \otimes V \rightarrow V, \ \ 
\theta (x \otimes y) = xy $, such that (\ref{trei}) and

 \begin{equation} \label{cinci}
xy  = yx  \ \ \ \ \forall x,y \in V
\end{equation} 
  
hold.

\bigskip

\begin{theorem}
Let  V   be a vector space spanned by   $a$ and  $b$, which are linearly independent.
Let                                                                                 $ \theta : V \otimes V \rightarrow V, \ \ 
\theta (x \otimes y) = xy $, 
be a linear map which satisfies (\ref{cinci}) and

 \begin{equation} \label{sase}
a^2 = b \ , \ \ \ b^2 = a \ .
\end{equation} 
                                                                               
Then:
$ ( V, \ \theta ) $
     is a Jordan algebra  
$ \iff $
$ ( V, \ \theta ) $
      is a non-unital commutative (associative) algebra.
\end{theorem}

{ \bf Proof. }
The indirect implication
               is obvious.

Let us prove the direct
implication.  
Since every vector in V can be expressed in terms of a  and  b, 
and using the commutativity of $ \theta $,
  we only need to check that

 \begin{equation} \label{sapte}
(ba)a = b a^2  \ , \ \ \ b^2 a = b (ba) \ ;
\end{equation}    
but these relations follow from (\ref{sase})
   and 
( \ref{trei}).
    Let us further observe that if
$ ba = \frac{1}{ \beta } a + \beta b \ , \ \  \beta^3 = -1 \ , $
then the relations ( \ref{sapte})   
are verified, otherwise  
$ ( V, \theta ) $
       is not a Jordan algebra. 
\qed

\bigskip
We will present the dual concept for a Jordan algebra below.

Let $ ( V, \ \theta ) $ be  a Jordan algebra, 
$ W \subset V \ot V \ot V \ot V $ the subspace generated by vectors
of the form $ y \ot x \ot x \ot x, \  x \ot y \ot x \ot x, \
x \ot x \ot y \ot x, \ x \ot x \ot x \ot y $, 
and $P: V ^{\ot 4} \rightarrow  W $ the projection associated to $W$.
A Jordan algebra satisfies the following relations:

$ \theta \circ  \tau  = \theta $

  $ \theta \circ ( \theta \ot I) \circ  ( \theta \ot I^{\ot 2} ) \circ P =
\theta \circ ( \theta \ot I) \circ  (I^{\ot 2} \ot \theta ) \circ P \ . $

\bigskip

In a dual manner, a Jordan coalgebra has a comultiplication
$ \eta : V \rightarrow V \ot V $ which satisfies:

$ \tau \circ \eta  = \eta $

 $ Q \circ ( \eta \ot I^{\ot 2}) \circ ( \eta \ot I) \circ  \eta =
Q \circ (I^{\ot 2} \ot \eta) \circ  ( \eta \ot I ) \circ \eta  \ , $

where $ Q: W \rightarrow  V ^{\ot 4} $ is the canonical inclusion
associated to W; in other words, the equality

$ ( \eta \ot I^{\ot 2}) \circ ( \eta \ot I) \circ  \eta =
(I^{\ot 2} \ot \eta) \circ  ( \eta \ot I ) \circ \eta  \ , $

should hold in W.

\bigskip

\cite{ ac} gives a theorem dual to
the Shirshov-Cohn theorem.
We will now state a theorem which is dual to the Theorem 2.1.

\bigskip

\begin{theorem}
Let  V   be a vector space of dimension two, and two linear maps
$ \eps, \ \zeta \ : V \rightarrow k $, which are linearly independent in $V^*$.
Let                                                                                 $ \eta : V \rightarrow V \ot V  \ $ 
be a linear map which satisfies
$ \  \tau \circ \eta = \eta \ , \ \ $
$(\eps \ot \eps ) \circ \eta = \zeta \ ,  \ \ 
(\zeta \ot \zeta ) \circ \eta = \eps $.

Then:
$ ( V, \ \eta ) $
     is a Jordan coalgebra  
$ \iff $
$ \eta $
      is cocommutative and coassociative.
\end{theorem}

\begin{remark}
The proof of the previous theorem follows by duality. Moreover, if $ {\  e, 
\ f \ } $ form a basis in $V$, then the conditions from hypothesis imply:

$ \eta (e) \ = \ \frac{1}{\beta} ( e \ot f + f \ot e ) \ + \ f \ot f \ $,

$ \eta (f) \ = \ {\beta} ( e \ot f + f \ot e ) \ + \ e \ot e \ $.

It is obvious that $ \eta $ is cocomutative.
The direct verification that  $ ( V, \ \eta ) $
     is a Jordan coalgebra is highly non-trivial.
Likewise, the 
 direct verification that $ \eta $ is coassociative is quite difficult.
\end{remark}

\bigskip

{\centering\section{THE QYBE AND ITS APPLICATIONS}}

For $ \  R: V \ot V \rightarrow V \ot V  $
a $ k$-linear map, let
$ {R^{12}}= R \ot I , \  {R^{23}}= I \ot R , \
{R^{13}}=(I\ot\tau )(R\ot I)(I\ot \tau ) $.

\begin{definition}
An invertible  $ k$-linear map  $ R : V \ot V \rightarrow V \ot V $
is called a Yang-Baxter
operator if it satisfies the  equation
\begin{equation}  \label{ybeq}
R^{12}  \circ  R^{23}  \circ  R^{12} = R^{23}  \circ  R^{12}  \circ  R^{23}
\end{equation}
An operator $R$ satisfies (\ref{ybeq}) if and only if
$R\circ \tau  $ satisfies
  the  QYBE 
(if and only if
$ \tau \circ R $ satisfies
the QYBE):
\begin{equation}   \label{ybeq2}
R^{12}  \circ  R^{13}  \circ  R^{23} = R^{23}  \circ  R^{13}  \circ  R^{12}
\end{equation}
\end{definition}

Examples of solutions to the QYBE from sets and Boolean algebras 
are described in \cite{HN}. Other examples related to differential geometry
are presented in  \cite{I, RI}. An exhaustive list of invertible solutions for (\ref{ybeq2}) in dimension 
2 is given in \cite{hi} and in the appendix of \cite{HlaSno:sol}.
Finding all Yang-Baxter operators in dimension greater than 2 is an 
unsolved problem. We will continue with more examples  below.

\bigskip

Let $A$ be a (unitary) associative $k$-algebra, and $ \alpha, \beta, \gamma \in k$. 
We define the
$k$-linear map:
$ \  R^{A}_{\alpha, \beta, \gamma}: A \ot A \rightarrow A \ot A, \ \ 
R^{A}_{\alpha, \beta, \gamma}( a \ot b) = \alpha ab \ot 1 + \beta 1 \ot ab -
\gamma a \ot b $.

\begin{theorem} (S. D\u asc\u alescu and F. F. Nichita, \cite{DasNic:yan})
\label{primat}
Let $A$ be an associative 
$k$-algebra with $ \dim A \ge 2$, and $ \alpha, \beta, \gamma \in k$. Then $ R^{A}_{\alpha, \beta, \gamma}$ is a Yang-Baxter operator if and only if one
of the following holds:

(i) $ \alpha = \gamma \ne 0, \ \ \beta \ne 0 $;

(ii) $ \beta = \gamma \ne 0, \ \ \alpha \ne 0 $;

(iii) $ \alpha = \beta = 0, \ \ \gamma \ne 0 $.

If so, we have $ ( R^{A}_{\alpha, \beta, \gamma})^{-1} = 
R^{A}_{\frac{1}{ \beta}, \frac{1}{\alpha}, \frac{1}{\gamma}} $ in cases (i) and
(ii), and $ ( R^{A}_{0, 0, \gamma})^{-1} = 
R^{A}_{0, 0, \frac{1}{\gamma}} $ in case (iii).
\end{theorem}

\begin{remark}

The Yang--Baxter equation plays an important role in knot theory. 
Turaev (\cite{T}) has described a general scheme to derive an invariant of 
oriented links from a Yang--Baxter operator.
In \cite{mn}, we considered the problem of applying Turaev's method to the
Yang--Baxter operators derived from algebra structures presented in
the above theorem. 
Turaev's procedure produced the 
Alexander polynomial of knots. 

\end{remark}

In dimension two, the operator from Theorem \ref{primat} (i) 
composed with the twist map, 
$ R^{A}_{\alpha, \beta, \alpha} \circ \tau$,
can be expressed as:
\begin{equation} \label{rmatcon2}
\begin{pmatrix}
1 & 0 & 0 & 0\\
0 & 1 & 0 & 0\\
0 & 1-q  & q & 0\\
\eta & 0 & 0 & -q
\end{pmatrix}
\end{equation}
where $ \eta \in \{ 0, \ 1 \} $, and $q \in k - \{ 0 \}$.
This form appears in the classifications of 
 \cite{hi, HlaSno:sol}.


\begin{definition}
A colored Yang-Baxter operator is defined as a function $ R
:X\times X \to \rend k {V\otimes V}, $ where $X$ is a set and $V$ is a
finite dimensional vector space over a field $k$. 
$R$
satisfies the two-parameter form of the QYBE if:
\begin{equation}\label{yb} 
R^{12}(u,v)R^{13}(u,w)R^{23}(v,w) = R^{23}(v,w)
R^{13}(u,w)R^{12}(u,v) \ \ \  \forall \ u,v,w\in X \ .
\end{equation} 
\end{definition}

\begin{theorem} (F. F. Nichita and D. Parashar, \cite{np}) \label{col}
Let $A$ be an associative 
$k$-algebra with $ \dim A \ge 2$, and
$ X \subset k $. Then,
for any two parameters $p,q\in k$, the function
$R:X\times X\to \rend k {A\otimes A}$ defined by
\begin{equation}\label{rsol} 
R(u,v)(a\otimes b) =p(u-v)1\otimes ab + q(u-v)ab\otimes 1 -(pu-qv)b\otimes a,
\end{equation}
satisfies the colored QYBE (\ref{yb}).
\end{theorem}

\begin{theorem} (F. F. Nichita and B. P. Popovici, \cite{nipo}) \label{top}
Let $A$ be an associative 
$k$-algebra with $ \dim A \ge 2$ and
 $q\in k$. Then the operator
\begin{equation}\label{slsol} 
S( \lambda )(a\otimes b) = (e^{\lambda} - 1)1\otimes ab 
+ q(e^{\lambda} - 1)ab\otimes 1 -(e^{\lambda}-q)b\otimes a
\end{equation}
satisfies the one-parameter form of the Yang-Baxter equation:
$$S^{12} (\lambda_{1} - \lambda_{2}) S^{13} (\lambda_{1} - \lambda_{3}) S^{23}(\lambda_{2} - \lambda_{3})=$$
\begin{equation}\label{onepara}
= S^{23} (\lambda_{2} - \lambda_{3}) S^{13} (\lambda_{1} - \lambda_{2})
S^{12} (\lambda_{1} - \lambda_{2}).
\end{equation}
\end{theorem}

\begin{remark} The operators from Theorems \ref{primat}, \ref{col}
and \ref{top} are related
via some algebraic operations, or the Baxterization procedure from \cite{defk}.
\end{remark}

\bigskip

Let $V$, $V'$, $V''$ be finite dimensional
vector spaces over $k$, and let $R: V\ot
V' \rightarrow V\ot V'$, $S: V\ot V'' \rightarrow V\ot V''$ and $T:
V'\ot V'' \rightarrow V'\ot V''$ be three linear maps.
The {\em Yang--Baxter
commutator} is a map $[R,S,T]: V\ot V'\ot V'' \rightarrow V\ot V'\ot
V''$ defined by \beq [R,S,T]:= R^{12} S^{13} T^{23} - T^{23} S^{13}
R^{12}. \eeq 

A system of linear
maps
$W: V\ot V\ \rightarrow V\ot V,\quad Z: V'\ot V'\ \rightarrow V'\ot
V',\quad X: V\ot V'\ \rightarrow V\ot V',$ is called a
$WXZ$--system if the
following conditions hold: \beq \label{ybsdoub} [W,W,W] = 0 \qquad
[Z,Z,Z] = 0 \qquad [W,X,X] = 0 \qquad [X,X,Z] = 0\eeq 
The above is one type of a constant
Yang--Baxter system that has been studied in \cite{np}, and
also shown to be closely related to entwining structures \cite{bn}.

\bigskip

\begin{theorem} (F. F. Nichita and D. Parashar, \cite{np})
Let $A$ be a $k$-algebra, and $ \lambda, \mu \in k$. The following is a 
$WXZ$--system:

$ W : A \ot A \rightarrow A \ot A, \ \ 
W(a \ot b)= ab \ot 1 + \lambda 1 \ot ab - b \ot a $,

$ Z : A \ot A \rightarrow A \ot A, \ \ 
Z(a \ot b)= \mu ab \ot 1 +  1 \ot ab - b \ot a $,

$ X : A \ot A \rightarrow A \ot A, \ \ 
X(a \ot b)= ab \ot 1 +  1 \ot ab - b \ot a $.

\end{theorem}

\begin{remark}
Let $R$ be a colored Yang-Baxter operator, i.e.
$ \ R^{12}(u,v)R^{13}(u,w)R^{23}(v,w) = R^{23}(v,w)
R^{13}(u,w)R^{12}(u,v)
 \ \  \forall \ u,v,w\in X$.

Then, if we let $ s, t \in X$, we obtain the following
$WXZ$--system:

$W= R(s, s) $,
$ \ X= R(s, t) $ and
$\ Z= R(t, t) $.

\end{remark}

\begin{definition}
A Lie superalgebra is a (nonassociative) $Z_2$-graded algebra, or superalgebra, 
over a field $k$ with the  Lie superbracket, satisfying the two conditions:
$$[x,y] = -(-1)^{|x||y|}[y,x$$
$$ (-1)^{|z||x|}[x,[y,z]]+(-1)^{|x||y|}[y,[z,x]]+(-1)^{|y||z|}[z,[x,y]]=0 $$
where $x$, $y$ and $z$ are pure in the $Z_2$-grading. Here, $|x|$ denotes the degree of $x$ (either 0 or 1). 
The degree of $[x,y]$ is the sum of degree of $x$ and $y$ modulo $2$.
\end{definition}

Let $ ( L , [,] )$ be a Lie superalgebra over $k$,
and  $ Z(L) = \{ z \in L : [z,x]=0 \ \ \forall \ x \in L \} $.
For $ z \in Z(L), \ \vert z \vert =0 $ and $ \alpha \in k $ we define:

$ \ \ \  { \phi }^L_{ \alpha} \ : \ L \ot L \ \ \longrightarrow \ \  L \ot L, \ \ 
x \ot y \mapsto \alpha [x,y] \ot z + (-1)^{ \vert x \vert \vert y \vert } y \ot x \ . $

Its inverse is:

$ \ \ \  {{ \phi }^L_{ \alpha}}^{-1} \ : \ L \ot L \ \ 
\longrightarrow \ \  L \ot L,  \ \ x \otimes y \mapsto \alpha z \otimes [x, y] + (-1)^{ \vert x \vert \vert y \vert 
} y \otimes x \ .$

\begin{theorem} (S. Majid, \cite{mj}) \label{majid}
Let  $ ( L , [,] )$ be a Lie superalgebra 
and
$ z \in Z(L), \vert z \vert = 0  $, and $ \alpha \in k $. Then:
$ \ \ \ \  { \phi }^L_{ \alpha} $ is a YB operator.
\end{theorem}

\begin{theorem}  (F. F. Nichita and B. P. Popovici, \cite{nipo}) \label{top2}
 Let  $ ( L , [,] )$ be a Lie superalgebra 
$ z \in Z(L), \vert z \vert = 0  $, 
$ X \subset k $,
and $ \alpha, \beta:X \rightarrow k $. 
Then,
$R:X\times X\to \rend k {L \otimes L}$ defined by
\begin{equation}\label{Lsol} 
R(u,v)(a\otimes b) = \alpha(u)[a,b]\otimes z + \beta (u) (-1)^{|a||b|} a\otimes b,
\end{equation}
satisfies the colored QYBE (\ref{yb}).
\end{theorem}

\begin{remark}
Let us consider the above data and apply it to Remark 3.9.
Then, if we let $ s, t \in X$, we obtain the following
$WXZ$--system:

$W(a\otimes b) = R(s,s)(a\otimes b)=  \ X(a\otimes b)= R(s,t)(a\otimes b)=
 \alpha(s)[a,b]\otimes z + \beta (s) (-1)^{|a||b|} a\otimes b, $ and

$Z(a\otimes b) = R(t,t)(a\otimes b)= 
\alpha(t)[a,b]\otimes z + \beta (t) (-1)^{|a||b|} a\otimes b$.

\end{remark}

\begin{remark}
The results presented in Theorems \ref{majid} and \ref{top2} hold for Lie algebras as well.
This is a consequence of the fact that these operators restricted to the
first component of a Lie superalgebra have the same properties.

\end{remark}

\bigskip
The constructions of this section were extended for 
 $(G,\theta)$-Lie algebra in \cite{nipo}.
For a  $(G,\theta)$-Lie algebra (see \cite{Kanak, nipo}) we have:
\begin{itemize}
\item $\langle L_a,L_b\rangle \subseteq L_{a+b} $
  \item $\theta$-braided (G-graded) antisymmetry: $\langle x,y\rangle = -\theta(a,b)\langle y,x  \rangle$ 
  \item $\theta$-braided (G-graded) Jacobi id: $\theta(c,a)\langle x, \langle y,z\rangle\rangle + \theta(b,c)\langle z, \langle x,y\rangle\rangle +\theta(a,b)\langle y, \langle z,x\rangle\rangle =0 $
  \item $\theta : G \times G \to C^*$ color function $\left \{ \begin{array}{c}\theta(a+b,c) = \theta(a,c)\theta(b,c)\\ \theta(a,b+c) = \theta(a,b)\theta(a,c)\\   \theta(a,b)\theta(b,a) =1\end{array}\right . $ 
\end{itemize}  

\bigskip

\begin{remark}

(i) In Theorem \ref{primat}, if we replace the associative algebra $A$, by a Jordan algebra $J$,
we obtain an operator which satisfies the braid equation if restricted to a subspace
$V= <a^2 \ot b \ot a, a\ot b \ot a^2 : a,b \in J >$ of $J^{\ot 3}$.

(ii) \cite{I, RI} present  construction of solutions for 
the Yang-Baxter equation  from Jordan triples and
from symmetric spaces.
Professor Dmitri Alekseevsky 
argued that these
are intimately related, and, in some cases, they might coincide. 

(iii) If we have in mind the results of \cite{bn, CMZ}, or the facts
presented in this section,
 the study of Jordan triples and the associated Yang-Baxter operators 
might lead to further constructions. The operators (\ref{slsol}) and (\ref{Lsol}) can be used to obtain $\theta$-dependent triple linear
products (see \cite{RI}, page 113); thus, they provide solutions for the
equation (2.2) of \cite{RI}.

(iv) Professor Takaaki Nomura pointed out that the Theorem 2.1
resembles the the Shirshov-Cohn Theorem. (The Shirshov-Cohn Theorem states that any Jordan algebra with two generators is special.) \cite{ac} 
presents a dual Shirshov-Cohn Theorem for Jordan coalgebras. 

(v) A fruitful observation (made at the 7-th Congress of Romanian Mathematicians) was that there are interesting connections between commutative Moufang loops and the Jordan identity. This is work in progress, and it is related to \cite{U}.

(vi) The Tits-Kantor-Koecher construction could be another way to relate the
Jordan algebras to the QYBE.
This can be done via the construction of
Yang-Baxter operators
from Lie algebras. Thus, the duality between Jordan algebras and Jordan
coalgebras is included in the self-duality of Yang-Baxter structures 
(see \cite{ns, nichita} for details about the self-duality of Yang-Baxter structures).
\end{remark}

\bigskip

\begin{center}

\end{center}

\bigskip

\bigskip

\bigskip

  Copyright © 2008 - 2012 by "Gheorghe Mihoc - Caius Iacob"
Institute of Mathematical Statistics and Applied Mathematics. All rights reserved.

\end{document}